\documentclass{tran-l}

\usepackage[active]{srcltx} % SRC Specials for DVI Searching

 \theoremstyle{definition}
 
 \theoremstyle{remark}
 
 \numberwithin{equation}{subsection}
\begin{document}

\title[Globally Solvable Vector Fields in Smooth Manifolds]
 {On The Existence of Globally Solvable Vector Fields in Smooth Manifolds}
\author{Jos\'e Ruidival dos Santos Filho }

\address{Departamento de Matem\'atica, Universidade Federal de S\~ao Carlos, S\~ao Carlos, 13565-905, SP, Brazil}
\email{santos@dm.ufscar.br}
\author{Joaquim Tavares}
\address{Departamento de Matem\'atica, Universidade Federal de Pernambuco, Recife, 50740-540, PE, Brazil}
\email{joaquim@dmat.ufpe.br}



%\date{00/00/000}

\dedicatory{}

\commby{}

%%% ----------------------------------------------------------------------

\begin{abstract}
Let $(\mathrm{M}, \omega_{0})$ be a connected paracompact smooth oriented manifold.  We establish a necessary and sufficient conditions on an involutive subbundle of $\mathrm{T\,M}$ such that $\mathrm{M}$ becomes simply connected.

\end{abstract}

%%% ----------------------------------------------------------------------
\maketitle
%%% ----------------------------------------------------------------------

\section{Introduction}

\medskip

\baselineskip 20pt

It is well known that there exist an obstruction to the existence of
$\mathrm{dim\,M}-k+1$ real linearly independent vector fields on an
manifold $\mathrm{M}$ in the $k$th cohomology group of
$\mathrm{M}$, the so called  $k$th Stiefel-Whitney class. We mean
that the $k$th Stiefel-Whitney class being nonzero implies that
there do not exist everywhere linearly independent vector fields. In
particular, the $\mathrm{dim\,M}$th Stiefel-Whitney class is the
obstruction to the existence of an everywhere nonzero vector field,
and the first Stiefel-Whitney class of a manifold is the obstruction
to orientability. Thus if one wishes to assume a hypotheses of
existence $k$ linearly independent vector fields in an orientable
manifold it certainly imposes the vanishing of such Stiefel-Whitney
classes. But it is not evident that additional hypothesis of
integrability and  global solvability of $k$ sections of
 $\Gamma({\textrm{T\,}\textrm{M}})$ the triviality of the all $l$th cohomology groups of
$\mathrm{M}$ for $l>\mathrm{dim\,}\textrm{M}-k+1$. This fact is encoded in a fundamental theorem of Hormander and Duistermaat
 (Theorem 6.4.2 [HD], pp 30) which characterizes global solvability of first order real differential operators in a $oriented$ smooth manifold.   But its proof is not trivial because it pass troughs a generalization of the theorem above. The generalization lies in a induction over the
 $convexity$ condition stated in $c)$ of Theorem 6.4.2 mentioned above. In the reading of proof one realize that the condition of convexity is needed only inside the orbits of the vector field
 as well as compactness is needed only for the intersection of the sublevels $\{x\in\mathrm{M}: u(x)\leq c\}$ with the orbits in $c)$ of Theorem 6.4.2.
 Taking this pointview we will extend the concept of convexity above by taking a  finitely generated Lie
 algebra $\mathfrak{a}\unlhd\Gamma(\textrm{T\,}\textrm{M})$ and  for a compact set $K\subset\textrm{M}$ we define ${\widehat{K}}$ to be the set obtained taking all smooth regular paths $\gamma$ with  with endpoints in $K$ and $\gamma^{\prime}\in\mathfrak{a}$.
We will show that in this context the analogous condition $c)$ of Theorem 6.4.2  is equivalent to one of the two condition stated below for the orbits
$\mathcal{O}_{\mathfrak{a}}(p)$ of $\mathfrak{a}$

(in the sense of Sussmann ([Su])).

\vskip0.5cm
(\texttt{C}) {\textit{ there exists at least $\mathrm{dim\,}\mathcal{O}_{\mathfrak{a}}(p)-1$ linearly independent
  sections of $\Gamma{(\mathrm{T\,}{\mathcal{O}_{\mathfrak{a}}(p)})}$
   which are   globally solvable for every $p\in\mathrm{M}$.}

   \vskip0.5cm
(\texttt{C}$^{\prime}$) {\textit{ there exist at least } $\mathrm{dim\,}\mathcal{O}_{\mathfrak{a}}(p)-1$
 {\textit{linearly independent sections  of }} $\Gamma{(\mathrm{T\,}{\mathcal{O}_{\mathfrak{a}}(p)})}$
  {\textit{ and }} $\mathcal{O}_{\mathfrak{a}}(p)${\textit{ is simply connected for every }$p\in\mathrm{M}$}.}
   \vskip0.5cm

  The condition (\texttt{C}$^{\prime}$) is a further restriction on the Stiefel-Whitney classes of $\textrm{T\,}\mathcal{O}_{\mathfrak{a}}(p))$ since we
  throughout this paper we will assume
$(\textrm{M}, \omega_{0})$ will be a paracompact manifold with positive orientation $\omega_{0}$. Then the condition (\texttt{C}$^{\prime}$) is entirely determined by the class of homotopy of the $CW$complex associated to $\text{M}$.
 Also we will assume the following hypothesis on the a finite set of generators $\{X_{1},...,X_{n}\}$ of $\mathfrak{a}\unrhd\Gamma(\mathrm{T\,}{\textrm{M}})$;
set $d_{l}$ to be the formal degree of $X_{l}$ and assume;
$$[X_{j}, X_{k}]=\sum_{d_{l}\leq d_{j}+d_{k}}  c_{jk}^{l}X_{l}$$
in the same sense as defined in the basic work of Nagel et al ([NSW]).  We also will denote by $\sim_{\mathfrak{a}}$ the equivalence relation of being in a same orbit $\mathcal{O}_{\mathfrak{a}}(p)$ and $\mathrm{H}_{\mathrm{d\,}}^{1}\mathcal{O}_{\mathfrak{a}}(p)$
the first De Rham cohomological group.
To avoid introduction of a metric in $\textrm{M}$ we will adopt the following definition found in page 111 of the reference [NSW].

{\textbf{Definitions.}}
Let $\mathfrak{a}$ a  Lie
 subalgebra of $\Gamma(\textrm{T\,}\textrm{M})$
 $finitely$ generated by $\{X_{1},...,X_{n}\}$ and assume that
  Denote by
 $C({\delta})$ the class of smooth paths $\gamma:[0,1]:\rightarrow\textrm{M}$ such that
 $$\gamma^{\prime}(s)=\sum^{n}_{l=1}c_{l}X_{l}(\gamma(s))$$ with
 $|c_{l}|\leq\delta^{d_{l}}$ and all $c_{l}$ constants.
 Define the pseudo-distance
 $$\rho(p,q)=\inf\{\delta>0: \exists\, \gamma\in C({\delta})\,\hbox{with}\,\gamma(0)=p, \gamma(1)=q\}$$
 and  denote by ${\widehat{K}_{\mathfrak{a}}}$ the  set of
  all paths $\gamma\in \cup_{\delta>0}C({\delta})$ with $\gamma(0), \gamma(1)\in K$.

{{We say that the triplet
 $(\textrm{M}, \mathfrak{a}, \omega_{0})$ is $\mathfrak
 {a}-$convex if for every compact subset $K\subset\textrm{M}$
 and $p\in\textrm{M}$ a pair of conditions holds;

{\textit{a) ${\widehat{K}_{\mathfrak{a}}}$ compact and}}

{\textit{b) there exists a set $\{X_{1},...,X_{n(p)-1}\}\in\Gamma{(\mathrm{T\,}{\mathcal{O}_{\mathfrak{a}}(p)})}$ of linearly independent sections such that
all its integral orbits are non
compact. }}

Our main result  is;

\textbf{Theorem.A.} Let $(\textrm{M}, \omega_{0})$ be a smooth oriented paracompact connected manifold and $\mathfrak{a}\unlhd\Gamma(\textrm{T\,}\textrm{M})$ be a finitely generated Lie algebra. Then  three statements are equivalent;

 A) $\textrm{M}$ is  $\mathfrak{a}-$convex,

 B) $\mathfrak{a}$ verifies   (\texttt{C}),

 C) $\mathfrak{a}$ verifies   (\texttt{C}$^{\prime}$).

\section{Proof of Theorem A.}
We start the proof with a lemmma;

{\textbf{Lemma.}\,\,}{\textit{Let $\mathfrak{a}\subset\Gamma(\mathrm{T}\mathrm{M})$ be a Lie subalgebra
   containing by $k-$linearly independent globally solvable vector fields. Then there exist a commutative subalgebra $\mathfrak{g}\unlhd\mathfrak{a}$ generated by $k-$linearly independent globally solvable vector fields.}

  \proof

 We perform induction on the number
  $k$.  For $k=1$ we just apply the equivalence between b) and f)
in the Theorem 6.4.2 in [DH]. Suppose that the Lemma is true for $k$ and let $\mathfrak{a}$ an Lie algebra containing
 $k+1$linearly independent set $\{X_{1},...,X_{n},X_{n+1}\}$ globally solvable vector fields.
  Then by induction hypothesis we may assume that $[X_{l},X_{k}]=0$ if $1\leq{l,k}\leq{n}$. Linear independence
  implies that $X_{n+1}=a_{1}X_{1}+\cdot\cdot\cdot+a_{n}X_{n}+Y$ with $Y\in\mathfrak{a}$ a non vanishing vector field. We define the a diffeomorphism   $\Phi:\text{M}\rightarrow\mathbf{R}^{n}\times{N}_{n}$ taking as the $n-$first coordinates of $\Phi$ the solutions $X_{k}t_{k}=1$ with its natural ordering and the trivialization comes again from the equivalence between a) and f) in Theorem 6.4.2 in [DH] just by a straightforward inductive argument over the number $k$. It follows that $\mathrm{D\,}\Phi(Y(p))\in\mathrm{T\,}_{p}({N}_{k})$ for every
  $p\in\textrm{M}$, where $\Pi^{-1}(\Pi(p))=\mathbf{R}^{n}\times \Pi(p)$ stands for  the  $n-k-$last coordinates of $\Phi$. Then by hypothesis for every
   orbit $\gamma$ of $X_{n+1}$ with endpoints $q_{0}, q_{1}\in K$ there exist another compact set $K^{\prime}\supset{K}$ such that $\gamma\subset{K}^{\prime}$. Then
   $\Pi(q_{0}), \Pi(q_{1})\in \Pi(K)$ are endpoints of  $\Pi\circ\gamma\subset\Pi({K}^{\prime})$ and $\mathrm{D}\Pi(\gamma^{\prime})=\mathrm{D}\Pi(Y(\gamma))\neq0$.
  If a integral orbit $\gamma$ of $X_{n+1}$ is not contained in a compact set $K$ then  the same should happens to  $\Pi(\gamma)$ with respect to the compact $\Pi(K)$. We now consider the orbits $\mathcal{O}_{\mathfrak{a}_{n}}(p)$ where $\mathfrak{a}_{n}\unlhd\mathfrak{a}$ is the Lie subalgebra generated by $\{X_{1},...,X_{n}\}$. All the orbits $\mathcal{O}_{\mathfrak{a}_{n}}(p)$ are diffeomorphic to $\mathbf{R}^{n}$.  We can solve $X_{n+1}u=1$ and consequently if
   $\gamma$ an integra
   l trajectory of $X_{n+1}$ parameterized  by $[0, \infty)$ with $\gamma(0)=p$ then

   $\lim_{t\rightarrow\infty}u(\gamma(t))=\infty$ and consequently $u$ will not be smooth in $\textrm{M}$ if $\Pi(\gamma)$ has compact closure.
  (Another argument: Denote by
    $\omega_{n}$ the positive orientation induced from $(\textrm{M},\omega_{0})$. Since
    $\omega_{n}(X_{1},...,X_{n})>0$ we can solve $X_{n+1}u=\omega_{n}(X_{1},...,X_{n})^{-1}$ entailing that
  $\mathrm{d\,}u\wedge\omega_{n}(X_{n+1},X_{1},...,X_{n})\equiv1$.
Let $K\subset{\textrm{M}}$ such that
     $K\cap\mathcal{O}_{\mathfrak{a}_{n}}(\gamma(t))$ is compact and
     $$1\,\mathrm{d\,}t=\int_{\Pi^{-1}
     (\Pi(\gamma(t)))}\chi_{K}\omega_{n}$$ for every $p\in\mathrm{M}$.
     If
   $\gamma$ an integral trajectory of $X_{n+1}$ parameterized  by $[0, \infty)$ with $\gamma(0)=p$ then
   $$
   \int_{\Pi(\gamma)}\bigg(\int_{\Pi^{-1}(\Pi(p))}
   \chi_{K}\omega_{n}\bigg)\mathrm{d\,}u
   =\int_{\Pi^{-1}(\Pi(\gamma))}
   \mathrm{d\,}u\wedge\chi_{K}\omega_{n}=\infty$$
   consequently if
   $\Pi(\gamma)$ if  compact closure because then the quantity on the left is finite}).  Thus any integral trajectory of $Y$ is unbounded and $Y$ verify the condition d) in Theorem 6.4.2 in [DH] is globally solvable in ${N}_{n}$ indeed and the equivalence between  f) and ) in Theorem 6.4.2 in [DH] applies to conclude the induction. $\square$

Now if B) is true we apply the Lemma to $\Gamma(\mathrm{T\,}\mathcal{O}^{\mathfrak{L}}(p))$ to find out that
is diffeomorphic
  to ${\mathbf{R}^{n(p)-1}}\times \gamma$ for some smooth regular connected curve $\gamma$ and with $n(p)=\mathrm{dim\,}\mathcal{O}^{\mathfrak{L}}(p)$.
    We consider  a non vanishing section of $\gamma^{\prime}\in\Gamma(\textrm{T}\gamma)$ which verifies $$\omega_{\mathcal{O}^{\mathfrak{L}}(p)}(\gamma^{\prime}(p), X_{1}(\gamma(p)),...,X_{n(p)-1}(p))>0$$ and define $$u(t)=\int^{t}_{ 0}\omega_{\mathcal{O}^{\mathfrak{L}}(s)}(\,\gamma^{\prime}(s), X_{1}(\gamma(s)),...,X_{n(p)-1}(\gamma(s))\,)\,\mathrm{d\,}s
   $$ is a strictly monotonous function of $t\in\mathbf{R}$ and consequently proper showing that $\gamma$ is diffeomorphic to an open interval and  the vector field $X_{\mathrm{dim\,}\mathcal{O}^{\mathfrak{L}}(p)}(p)=
   \gamma^{\prime}(p)$ is indeed globally solvable in $\mathrm{M}$ by d) in Theorem 6.4.2 in [DH]. Then we apply again the Lemma to find out that ${\mathcal{O}_{\mathfrak{a}}(p)}
 \simeq\textbf{R}^{\textrm{dim\,}{\mathcal{O}_{\mathfrak{a}}(p)}}$.
Observe that ${\mathcal{O}_{\mathfrak{a}}(p)}\cap{\widehat{K}_{\mathfrak{a}}}=
 \widehat{{\mathcal{O}_{\mathfrak{a}}(p)}\cap{{K}_{\mathfrak{a}}}}$  comes from a result of Sussmann ([Su]) and consequently compactness of
$\widehat{K}_{\mathfrak{a}}$ is equivalent to compactness of its intersection with an arbitrary orbit ${\mathcal{O}_{\mathfrak{a}}(p)}$.
 On the other hand the condition $(\texttt{C})$
 implies that ${\mathcal{O}_{\mathfrak{a}}(p)}
 \simeq\textbf{R}^{\textrm{dim\,}{\mathcal{O}_{\mathfrak{a}}(p)}}$ with $\widehat{K}_{\mathfrak{a}}\cap{\mathcal{O}_{\mathfrak{a}}(p)}$ corresponding to the  convex envelope of $K\cap{\mathcal{O}_{\mathfrak{a}}(p)}$ in
 $\textbf{R}^{\textrm{dim\,}{\mathcal{O}_{\mathfrak{a}}(p)}}$. Since  the Lie algebra $\mathfrak{a}$ is finitely generated and ${\mathcal{O}_{\mathfrak{a}}(p)}
 \simeq\textbf{R}^{\textrm{dim\,}{\mathcal{O}_{\mathfrak{a}}(p)}}$  the metric $\rho$  is well defined by $\mathfrak{a}$ in ${\mathcal{O}_{\mathfrak{a}}(p)}$. Denote by $\mathrm{diam\,}$ the set function induced in ${\mathcal{O}_{\mathfrak{a}}(p)}$
 by the Euclidean  metric in $\textbf{R}^{\textrm{dim\,}{\mathcal{O}_{\mathfrak{a}}(p)}}$ with the inherited orientation from $\mathrm{M}$ (see Theorem 1. \& 3, page 110 in [NSW]).  It follows that if $p,q\in K\cap\mathcal{O}_{\mathfrak{a}}(p)$  then
 $$C_{1}\, \textrm{diam}(K\cap\mathcal{O}_{\mathfrak{a}}(p))\leq\rho(q, p)\leq{C_2}\,\textrm{diam\,}^{1/\max d_{j}}(K\cap\mathcal{O}_{\mathfrak{a}}(p))
$$
 from Proposition 1.1 page 107 in [NSW]. Consequently $$\widehat{K}_{\mathfrak{a}}\cap\mathcal{O}_{\mathfrak{a}}(p)
 \subset
 \{q\in\mathcal{O}_{\mathfrak{a}}(p):\rho(q, p)\leq
C_{3}\textrm{diam\,}^{1/\max d_{j}}(K\cap\mathcal{O}_{\mathfrak{a}}(p))\}$$
and
 $\{\rho(\cdot, p)\leq c\}\cap{\mathcal{O}_{\mathfrak{a}}}(p)$ is always compact for every $c\in\mathbf{R}$. Moreover if $(x_{1},...,x_{n(p)})$ are the coordinates of $\mathbf{R}^{\textrm{dim\,}{\mathcal{O}_{\mathfrak{a}}(p)}}$ then the inverse image of  $(\partial/\partial x_{1},...,\partial/\partial x_{n(p)})$ in $\Gamma{(\mathrm{T\,}{\mathcal{O}_{\mathfrak{a}}(p)})}$ will verify the property $b)$ for $\mathfrak{a}-$convexity, completing the first part of the proof.

Now if A) is true, since $\mathcal{O}_{\mathfrak{a}}(p)$ is a manifold and locally we may assume that $V\simeq\mathbf{R}^{\textrm{dim}\mathcal{O}_{\mathfrak{a}}(p)}$
for an open $convex$ neighborhood $V\subset{\textrm{M}}$ of $p\in\textrm{M}$. Then $(V,\mathfrak{a}, \omega_{0})$ verifies
$(\texttt{C})$ and
$\{\rho(\cdot, p)\leq c\}\cap{\,\mathcal{O}_{\mathfrak{a}}}(p)$ remains compact for small $c\in\mathbf{R}$.
 On the other hand by besides the existence of global linearly independent set  $\{X_{1},...,X_{n(p)-1}\}$ in $\Gamma(\mathrm{T\,}\mathcal{O}_{\mathfrak{a}}(p))$ by the property $b)$ of $a-$convexity, the pseudo-distances defined by $\rho$ but taking only smooth paths
 tangent to the linear span of the set $\{X_{1},...,X_{n}\}$ generating $\mathfrak{a}$ is equivalent to the former one by
 Theorem 3. \&4. Equivalent pseudo-distances, page 111 in [NSW].
 As a consequence the set $\{X_{1},...,X_{n(p)-1}\}$ will be a set of
 globally solvable vector fields of $\Gamma(\mathrm{T\,}\mathcal{O}_{\mathfrak{a}}(p))$ by the property $a)$ of $\mathfrak{a}-$convexity and the condition {(\texttt C)} is verified and A) implies B). That B) implies C) is just application of the Lemma together the equivalence  $f)$ in in Theorem 6.4.2 in [DH].
 Then it is left to prove that C) implies B). Since
 $\mathfrak{a}$ is finitely generated we may select a set
 $\{X_{1},...,X_{n}\}\subset \Gamma(\textrm{T\,}\textrm{M})$ generating $\mathfrak{a}$ and define the second order operator
 $H=X^{2}_{1}+\cdot\cdot\cdot+X^{2}_{n}$ which is well defined in every orbit ${\mathcal{O}_{\mathfrak{a}}}(p)$ where it is hypoelliptic by a result of Hormander[Ho]) . It is a consequence of the Bony's maximum principle ([Bo]) that  a twice differentiable function $u$ verifying  $Hu\geq 0$ have upperlevel  sets
 $\{u\geq c\}\cap{\mathcal{O}_{\mathfrak{a}}}(p)$ is relatively non compact in $\mathcal{O}_{\mathfrak{a}}$. On the other hand if we assume also that the sublevel  sets
 $\{u\leq c\}\cap{\mathcal{O}_{\mathfrak{a}}}(p)$ are compact
 and the $2$th Stiefel-Whitney class of
 $\textrm{T}\,{\mathcal{O}_{\mathfrak{a}}}(p)$  is trivial  there exist $n(p)-1$ linearly independent sections $\{X_{1},...,X_{n(p)-1}\}\subset
 \Gamma(\textrm{T\,}{\mathcal{O}_{\mathfrak{a}}}(p))$.
 Since the orbits ${\mathcal{O}_{\mathfrak{a}}}(p)$ are positively oriented with the inherited
  $n(p)-$differential form $\omega_{{\mathcal{O}_{\mathfrak{a}}(p)}}$, the $1-$differential form
 $$\omega_{\mathcal{O}_{\mathfrak{a}}(p))}
 (\,\cdot\,,X_{1},...,X_{n(p)-1})$$
 is nonvanishing. But it  vanishes in the tangent space of the orbit generated by the the linearly independent set of  vector fields $\{X_{1},...,X_{n(p)-1}\}$ which dimension is $n(p)-1$ or $n(p)$, the latter incompatible with the positivity of the orientation in $\mathcal{O}_{\mathfrak{a}}(p)$.
 If the first De Rham cohomology group of
 $\mathcal{O}_{\mathfrak{a}}(p)$ is trivial then locally we can find smooth function $u$ with  $\mathrm{d\,}u=\omega_{\mathcal{O}_{\mathfrak{a}}(p)}
 (\,\cdot\,,X_{1},...,X_{n(p)-1})$ which is constant in the components of  orbits generated by $\{X_{1},...,X_{n(p)-1}\}$. Then for any choice of $X_{n(p)}\in \Gamma(\textrm{T\,}{\mathcal{O}_{\mathfrak{a}}}(p))$ linearly independent from
 $\{X_{1},...,X_{n(p)-1}\}$ such that
 $\omega_{\mathcal{O}_{\mathfrak{a}}(p)}
 (X_{n(p)},X_{1},...,X_{n(p)-1})>0$
 will verify
 $X^{2}_{n(p)} \mathrm{e}^{-\kappa\,u^{2}}(p)> 0$ for large positive $\kappa(K)$ and all $p\in K$, a compact set. The paracompactness of $\textrm{M}$ allows one to find sequence of compacts $K_{n}$ such that $\textrm{int\,}K_{n}\subset K_{n+1} $ and smooth $\chi_{K_{n}}$ with $\chi_{K_{n}}(q)=1$ if
 $q\in K_{n}$, $\chi_{K_{n}}(q)=0$ if $q\in \text{M}\setminus{K_{n+1}}$ and $0\leq \chi_{K_{n}}(q)\leq 1$ for all $q\in\textrm{M}$. Denote by $\kappa(K_{n})$ a constant
 verifying
 $$\kappa(K_{n})\min_{p\in{K}_{n}}\omega_{\mathcal{O}_{\mathfrak{a}}(p)}
 (X_{n(p)},X_{1},...,X_{n(p)-1})\geq \sup_{p \in{K}_{n-1}}
 |X^{2}_{n(p)}\chi_{K_{n}}u|.$$
 Then by a suitable choice of constants $\kappa(K_{n})$ we find out that the function $$u_{0}=\sum^{\infty}_{n=1}\chi_{K_{n}} \mathrm{e}^{-\kappa(K_{n})\,u^{2}}$$
 is smooth and $X^{2}_{n(p)}u_{0}>0$ for all $p\in\textrm{M}$. If the
integral
trajectory $\gamma$ of $X_{n(p)}$ starting at $p\in \mathcal{O}_{\mathfrak{a}}(p))$ remains in a compact set $K$ then at some point $q$ of this compact $$\omega_{\mathcal{O}_{\mathfrak{a}}(p)}
(X_{n(p)}(q),X_{1}(q),...,X_{n(p)-1}(q))=0$$ which is a contradicts the positivity of the orientation $\omega_{\mathcal{O}_{\mathfrak{a}}(p))}(q)$.
It follows that the condition $c)$ in
 in Theorem 6.4.2 in [DH] is verified for  $X_{n(p)}$ and it is globally integrable and we apply f) in the same theorem to write
 $\mathcal{O}_{\mathfrak{a}}(p))=N_{n(p)-1}\times\textbf{R}$.
 Since the submanifolds $N_{n(p)-1}\times {t}$ inherits same
 properties of  $\mathcal{O}_{\mathfrak{a}}(p)$  we apply induction to conclude that $\mathcal{O}_{\mathfrak{a}}(p)\simeq\mathbf{R}^{n(p)}$,
 finishing the proof. $\square$

Before finish this paper we must remark the the extraordinary semblance between this form of
presenting global solvability for a real vector field in a manifold with the
parallelization Theorem A  in a work of
 Greene \& Shiohama ([GS]), but in the absence of a Riemannian structure or
the necessity of emptiness of the critical set of $u$.
The Theorem B in [GS] says that when the manifold is not simply
connected the obstruction to global solvability will be located
in the singular set of the convex function pointing to further
investigation of this phenomena. Also one also may apply the generalized Tietze-Nakajima theorem in [KB] together the Lemma above to show that  the manifold $\textrm{M}$ has as projection
of the first $n-$coordinates of $\Phi$ a convex subset of $\textbf{R}^{n}$ when $\Phi$ is proper.

\end{document}